\documentclass[reqno]{amsart}

\usepackage{amsmath}
\usepackage{amssymb}
\usepackage{amsfonts}
\usepackage{mathrsfs}
\usepackage{bbm}

\usepackage{todonotes}

\usepackage[tmargin=1.5in,bmargin=1.5in,lmargin=1.3in,rmargin=1.3in]{geometry}
\usepackage{microtype}

\newtheorem{thm}{Theorem}

\newtheorem{lem}[thm]{Lemma}

\theoremstyle{definition}

\theoremstyle{remark}

\numberwithin{equation}{section}

\usepackage[
bookmarks, colorlinks, breaklinks, pdftitle={The critical surface of quasi-transitive inhomogeneous percolation models},pdfauthor={Thomas Beekenkamp and Tim Hulshof}]{hyperref}  
\hypersetup{linkcolor=black,citecolor=black,filecolor=black,urlcolor=black} 

\title{Sharpness for inhomogeneous percolation on quasi-transitive graphs}

\author{Thomas Beekenkamp}
\address{Mathematisches Institut, Ludwig-Maximilians-Universit\"at M\"unchen, Theresienstra\ss{}e 39, 80333 M\"unchen, Germany}
\email{Thomas.Beekenkamp@math.lmu.de}

\author{Tim Hulshof}
\address{Department of Mathematics and Computer Science, Eindhoven University of Technology, P.O. Box 513, 5600 MB Eindhoven, The Netherlands.}
\email{w.j.t.hulshof@tue.nl}

\newcommand\nlongleftrightarrow{\mathrel{\,\,\,\not\!\!\!\longleftrightarrow}}

\newcommand{\Gcol}{G_{\mathrm{col}}}
\newcommand{\Ppq}{\mathbb{P}_{\mathbf{p},q}}
\newcommand{\Ccal}{\mathscr{C}}
\newcommand{\ssubset}{\subset\!\subset}
\newcommand{\sss}{\scriptscriptstyle}

\begin{document}
\begin{abstract}
In this note we study the phase transition for percolation on quasi-transitive graphs with quasi-transitively inhomogeneous edge-retention probabilities. A quasi-transitive graph is an infinite graph with finitely many different ``types'' of edges and vertices. We prove that the transition is sharp almost everywhere, i.e., that in the subcritical regime the expected cluster size is finite, and that in the subcritical regime the probability of the one-arm event decays exponentially. Our proof extends the proof of sharpness of the phase transition for homogeneous percolation on vertex-transitive graphs by Duminil-Copin and Tassion \cite{Duminil20152}, and the result generalizes previous results of Antunovi\'c and Veseli\'c \cite{Antunovic2008} and Menshikov \cite{Menshikov1986}.
\end{abstract}
\maketitle
\section{Introduction and main result}
Roughly speaking, this note is about independent percolation models on graphs with enough symmetry that there are only a finite number of ``types'' of vertices and edges, and where each type of edge receives its own percolation parameter. In this setting, there are often many different combinations of the parameters that give rise to critical behavior, and when embedded in Euclidean space, we can view these parameter combinations as points on a surface. Our aim is to study properties of this surface and the regions of the parameter space it separates. 

To make the above notions precise, and to state our main result, we will need to start with some notation and definitions:

 We define an \emph{edge-coloured graph} $\Gcol$ with $N \in \mathbb{N}$ colours as a (countable) set $V$ and a collection of $N$ mutually disjoint sets $E_1, \dots, E_N$ of subsets of pairs of elements from $V$. We write $\Gcol =(V, E_1, \dots, E_N)$. We call $V$ the \emph{vertex set} and the union $E:= \cup_{i=1}^N E_i$ the \emph{edge set} of $\Gcol$.
In this note we consider only infinite, connected, locally finite, edge-coloured graphs $\Gcol$, i.e., satisfying $|V| = \infty$, (where $| \cdot |$ denotes the cardinality of a set) such that there exists no element $v \in V$ that is contained in an infinite number of edges, and where between each two vertices there exists a path of edges in the edge set $E$.
Moreover, we only consider quasi-transitive edge-coloured graphs:

Let $\Gcol = (V, E_1, \dots, E_N)$ be an infinite, locally finite, connected coloured graph. We call a bijective map $f : V \to V$ a \emph{coloured graph automorphism} if $\{u,v\} \in E_i$ if and only if $\{f(u),f(v)\} \in E_i$ for all $1 \le i \le N$.
We say $\Gcol$ is \emph{quasi transitive} if we can partition $V$ into $k \in \mathbb{N}$ disjoint subsets $V = \dot \cup_{i=1}^k V_i$ such that for any $1 \le j \le k$ and any $u,v \in V_j$ there exists a coloured graph automorphism $f$ such that $f(v)=u$. Note that the monochromatic case $N=1$ corresponds to the usual notion of a quasi-transitive graph (see e.g. \cite{Haggstrom}) and if moreover $k$ can be chosen equal to $1$, it corresponds to the usual notion of a vertex-transitive graph.

We use edge-coloured graphs to define a quasi-transitive inhomogeneous bond percolation model, where each colour edge has a different edge retention probability:

Given $\Gcol = (V, E_1,\dots,E_N)$, $\mathbf{p} := (p_1,\dots, p_{N-1}) \in [0,1]^{N-1}$ and $ q := p_N \in [0,1]$, the \emph{inhomogeneous percolation measure} $\Ppq$ is the measure on $\{0,1\}^E$ where each edge of $E_i$ \emph{retained} with probability $p_i$, and \emph{removed} with probability $1-p_i$.
More precisely, we consider the probability space $(\Omega, \mathcal{F}, \Ppq)$ where $\Omega = \prod_{e \in E} \{0,1\}$ is the sample space and $\mathcal{F}$ is the $\sigma$-algebra of $\Omega$ generated by the cylinder events.
For $1 \le i \le N$ and $e \in E_i$, let $\mu_e^{\sss (i)}$ denote a Bernoulli measure on $\{0,1\}$ with parameter $p_i$, i.e.,
\[
	\mu_{e}^{\sss (i)}(\omega(e) = 1) = p_i \qquad \text{ and } \qquad \mu_e^{\sss (i)}(\omega(e)=0) = 1-p_i.
\]
We define $\Ppq$ associated with $\Gcol$ to be the product measure
\[
	\Ppq := \prod_{i=1}^N \prod_{e \in E_i} \mu_e^{\sss (i)}.
\]
We refer to the ``degenerate'' case of an inhomogeneous percolation measure on a graph where all edges are the same colour (i.e., $N=1$) as \emph{homogeneous percolation.} (Note that this is more commonly known as \emph{Bernoulli bond percolation,} or, considering that this is by far the most studied class of models in percolation theory, simply as ``percolation''.)

In the definition of $\Ppq$ we have singled out the $N$-th parameter and called it $q$. We did this because we are interested in the transitions associated with the cluster from finite to infinite. Since the probability of an infinite cluster and the expected size of a cluster cannot decrease when we increase any parameter, we can view these transitional points as forming a surface in $[0,1]^N$ that is single-valued in any of its parameters, so we can use $q$ to parametrize the surface (but there is nothing ``special'' about $q$).

We write $\Ccal(x)$ for the connected component of $x \in V$ in a percolation configuration. We define the \emph{percolation function}
\[
	\theta(\mathbf{p},q) := \sup_{x \in V} \Ppq(|\Ccal(x)| =\infty)
\]
and the \emph{susceptibility:}
\[
\chi(\mathbf{p},q) :=\sup_{x\in V}\mathbb{E}_{\mathbf{p},q}[|\mathscr{C}(x)|].
\]
Associated with $\theta(\mathbf{p},q)$ is the \emph{$\theta$-critical surface}
\[
q_\theta(\mathbf{p}) :=\inf\left\{q\in [0,1] \::\: \theta(\mathbf{p},q) >0\right\}
\]
and the associated with $\chi(\mathbf{p},q)$ is the \emph{$\chi$-critical surface:}
\[
q_\chi(\mathbf{p}) :=\inf\big\{q\in[0,1] \::\:\chi(\mathbf{p},q)=\infty\big\}.
\]
Note that $\chi(\mathbf{p},q) =\infty$ whenever $\theta(\mathbf{p},q) >0$, so $q_\chi(\mathbf{p}) \le q_\theta(\mathbf{p})$. 

The classical \emph{sharpness of the phase-transition} question is whether it holds that $q_\chi(\mathbf{p}) = q_\theta(\mathbf{p})$. In general this is not true, as it is easy to find counter-examples. But these counter-examples are often considered pathological, and thus, uninteresting. The question of sharpness becomes interesting, however, if we restrict ourselves to more ``reasonable'' graphs. In the current paper, we therefore consider only quasi-transitive coloured graphs with quasi-transitive inhomogeneous percolation measures. 

We will show sharpness of the phase transition for inhomogeneous percolation models for almost every $\mathbf{p}$. That is, we show that the two critical surfaces $q_\theta(\mathbf{p})$ and $q_\chi(\mathbf{p})$ are equal almost everywhere. The proof of this relies on a different characterisation of the critical surface, which we will introduce first. 

We write $d(u,v)$ for the graph distance between two vertices $u,v\in V$, i.e., for the number of edges in the shortest path from $u$ to $v$ in $\Gcol$. We write $\Lambda_k^u$ for the ball of radius $k$ centered at $u$, and $\partial \Lambda_u^k$ for its vertex boundary, i.e.,
\[
\Lambda_k^u:=\{v\in V\::\:d(u,v)\leq k\} \qquad \text{ and } \qquad \partial \Lambda_{k}^u:=\Lambda_k^u \setminus \Lambda_{k-1}^u.
\]

This brings us to our main result:
\begin{thm}\label{thm}
Let $G_{\mathrm{col}}= (V,E_1,\dots, E_N)$ be an infinite, locally finite, connected, quasi-transitive coloured graph with $N$ colours. Then:
\begin{enumerate}
	\item for almost all $\mathbf{p}\in [0,1]^{N-1}$ the transition is sharp, i.e., $q_\theta(\mathbf{p}) = q_\chi(\mathbf{p})$, and
	\item for almost all $\mathbf{p}\in [0,1]^{N-1}$, if $q<q_\theta(\mathbf{p})$, then there exists a constant $c>0$ such that for all $x\in V$ and for all $k\in \mathbb{N}$
	\begin{equation}\label{eq:thma}
	\mathbb{P}_{\mathbf{p},q}\left(x\longleftrightarrow \partial\Lambda_k^x\right)\leq e^{-ck}.
	\end{equation}
\end{enumerate}
\end{thm}

By `for almost all' we mean that the set of $\mathbf{p}\in [0,1]^{N-1}$ for which the transition is not sharp has Lebesgue-measure 0. We expect Theorem \ref{thm} to hold for all $\mathbf{p}\in [0,1]^{N-1}$, however proving this will be difficult, in light of the results by Balister, Bollob\'as and Riordan in \cite{balister2014essential}. There the authors have an inhomogeneous percolation model with two parameters. They show that it is not clear, except in a few special cases, if an increase of one parameter will decrease the critical value of the other parameter. The case where the critical value does not decrease corresponds to a discontinuity in $q_\theta$. We will show below that the set of these discontinuities has measure 0.

Theorem \ref{thm} is in various ways an extension of what is known in the literature. In particular, Menshikov \cite{Menshikov1986} proved sharpness of the phase transition for homogeneous percolation on vertex-transitive graphs where the boundaries of metric balls rooted at fixed vertices grow subexponentially. As Menshikov remarks, his proof can be extended (with some effort) to the inhomogeneous setting described in this paper, but the subexponential growth condition is essential to his proof. The subexponential growth condition excludes a large class of graphs from his result, for instance, all graphs that are isometric to a $d$-regular tree, or to a hyperbolic space, or to a non-amenable (Cayley) graph. 
Around the same time as Menshikov, on the other side of the iron curtain, Aizenman and Barsky \cite{AizenmanBarsky1987}, using a different method, proved that the phase transition is sharp for various lattice models, again requiring something like the subexponential growth condition. Antunovi\'c and Veseli\'c \cite{Antunovic2008} extended the latter proof to homogeneous percolation models on all quasi-transitive graphs, including those with boundary growth that is not sub-exponential. Most recently, Duminil-Copin and Tassion \cite{Duminil20152} introduced an entirely different method for proving sharpness of the phase transition for (among others) homogeneous percolation on vertex-transtive graphs. (They also give a more streamlined version of the proof for nearest-neighbor percolation on $\mathbb{Z}^d$ \cite{duminil2016new}.)
We show how to extend the proof of Duminil-Copin and Tassion from uniqueness at the critical point of homogeneous percolation to uniqueness of the critical surface of inhomogeneous percolation.

\section{Outline of the proof}
The core of the proof is to define an alternative critical surface, $q_\psi(\mathbf{p})$. We then use this alternative characterization to show that all three critical surface coincide. To give the alternative definition we first need some notation. For a set $S\subset V$ we write $\Delta S$ for its edge boundary, i.e., $\{x,y\}\in \Delta S$ if $x\in S$ and $y\not\in S$ and $\{x,y\}\in E$. Furthermore for $S\subset V$ we write $x\stackrel{S}{\longleftrightarrow}y$ for the event that there is a path of open edges that have both endpoints in $S$ from $x$ to $y$. We write $S \ssubset V$ to denote that $S \subset V$ and $|S| < \infty$.
We define for $x\in V$ and for $S \ssubset V$ such that $S \owns x$ (recalling that $\mathbf{p} = (p_1, \dots, p_{N-1})$ and $q = p_N$),
\begin{equation}\label{e:psidef}
\psi_{\mathbf{p},q}(x,S):=\sum_{i=1}^{N} p_i\sum_{\substack{\{y,z\}\in \\ \Delta S\cap E_i}}\mathbb{P}_{\mathbf{p},q}\left(x\stackrel{S}{\longleftrightarrow} y\right).
\end{equation}
(Note that $\psi_{\mathbf{p},q}(x,S)$ is the analogue of $\phi_p(S)$ defined in \cite{Duminil20152} for homogeneous percolation on transitive graphs.) We define the \emph{$\psi$-critical surface}
\begin{align}\label{eq:qt}
q_\psi(\mathbf{p}):=\sup\Big\{q\in [0,1]\::\:\forall x\in V \:\exists S\ssubset V\text{ with } x\in S, \text{ and }\psi_{\mathbf{p},q}(x,S)<1\Big\},
\end{align}
or equivalently,
\begin{equation}\label{eq:qt2}
q_\psi(\mathbf{p})= \inf_{x\in V} \sup_{\substack{S\ssubset V: \\ x\in S}} \sup\big\{q\in[0,1]\: :\: \psi_{\mathbf{p},q}(x,S)<1\big\}.
\end{equation}

The aim is now to show that $q_\psi=q_\theta=q_\chi$. To start we prove that from $q<q_\psi(\mathbf{p})$ it follows that $\chi(\mathbf{p},q) < \infty.$ From that we can conclude $q_\psi(\mathbf{p}) \leq q_\chi(\mathbf{p})\leq q_\theta(\mathbf{p})$. We show that $\chi(\mathbf{p},q)$ is indeed finite when $q<q_\psi(\mathbf{p})$ in Section \ref{sec:fin} below.

For the supercritical part of the proof we consider values of $\mathbf{p}\in [0,1]^{N-1}$ at which $q_\psi$ is continuous in the following sense: we consider values of $\mathbf{p}$ in the set
\[
\mathcal{A}:=\left\{\mathbf{p}\in [0,1]^{N-1} \::\: \lambda \mapsto q_\psi(\mathbf{p+\lambda \mathbf{1}})\text{ is left-continuous at }\lambda=0\right\},
\]
where $\mathbf{1}$ denotes the vector of length $N-1$ with all entries equal to $1$. If for some $x\in V$ and all $q>q_\psi(\mathbf{p})$,
\begin{equation}\label{eq:super}
\mathbb{P}_{\mathbf{p},q}(|\Ccal(x)| = \infty)>0,
\end{equation}
then it follows that $q_\psi(\mathbf{p}) \ge q_\theta(\mathbf{p})$. We prove that \eqref{eq:super} holds when $q>q_\psi(\mathbf{p})$ using the continuity of $q_\psi(\mathbf{p})$ in Section \ref{sec:sup} below. It is easy to show that $A^c$ has measure $0$: It holds that $q_\psi(\mathbf{p}+\lambda \mathbf{1})$ is decreasing in $\lambda$. It follows that there are at most countably many values of $\lambda$ for which $q_\psi(\mathbf{p}+\lambda \mathbf{1})$ is discontinuous. If we want to integrate $\mathbbm{1}_{\mathcal{A}^c}$ with respect to $p_1,\dots,p_{N-1}$ we can do a change of variables and first integrate with respect to $\lambda$ and then with respect to $p_1,\dots,p_{N-1}$. This gives
\[
\int \mathbbm{1}_{\mathcal{A}^c}\mathrm{d}\mathbf{p}=\idotsint C \mathbbm{1}_{\{q_\psi(\mathbf{p}+\lambda\mathbf{1}) \textrm{ is not left-continuous at } 0\}}\mathrm{d}\lambda\mathrm{d}p_1\cdots\mathrm{d}p_{N-1}=0.
\]

In Section \ref{sec:exp} we prove that \eqref{eq:thma} holds for $q<q_\psi(\mathbf{p})$. Combining this with Theorem \ref{thm}(a) gives Theorem \ref{thm}(b).

\section{Subcritical Case: Finite Susceptibility}\label{sec:fin}
In this section we prove finite susceptibility below $q_\psi(\mathbf{p})$ for any $\mathbf{p}\in [0,1]^{N-1}$. Let $q<q_\psi(\mathbf{p})$. Therefore, for any $x\in V$ there exists a finite set $S_x\subset V$ with $x\in S_x$ and $\psi_{\mathbf{p},q}(x,S_x)<1$. For a finite set $A\ssubset V$ define
\begin{equation}
\chi(A)=\max_{u\in V}\sum_{v\in A}\mathbb{P}_{\mathbf{p},q}\left(u\stackrel{A}{\longleftrightarrow}v\right).
\end{equation}
Let $x \in V$ and let $A\ssubset V$ be such that $A\setminus S_x\neq \varnothing$. Now suppose that the event $x\stackrel{A}{\longleftrightarrow}u$ holds for some $u\in A\setminus S_x$. Then there exists an open edge $e=\{y,z\}$ on the boundary $\Delta S_x$ such that $x\stackrel{S_x}{\longleftrightarrow} y$ and $z\stackrel{A}{\longleftrightarrow} u$ using disjoint paths. We can use this observation together with the BK inequality \cite{berg1985} to bound the probability of $x\stackrel{A}{\longleftrightarrow}u$:
\begin{equation}
\mathbb{P}_{\mathbf{p},q}\left(x\stackrel{A}{\longleftrightarrow}u\right)\leq \sum_{i=1}^N\sum_{\substack{\{y,z\}\in \\ \Delta S_x\cap E_i}}p_i\mathbb{P}_{\mathbf{p},q}\left(x\stackrel{S_x}{\longleftrightarrow} y\right)\mathbb{P}_{\mathbf{p},q}\left(z\stackrel{A}{\longleftrightarrow} u\right).
\end{equation}
Summing the above inequality over all $u\in A\setminus S_x$ gives:
\[
\begin{split}
\sum_{u\in A\setminus S_x}\mathbb{P}_{\mathbf{p},q}\left(x\stackrel{A}{\longleftrightarrow}u\right)&\leq \sum_{u\in A\setminus S_x}\sum_{i=1}^N\sum_{\substack{\{y,z\}\in \\ \Delta S_x\cap E_i}}p_i\mathbb{P}_{\mathbf{p},q}\left(x\stackrel{S_x}{\longleftrightarrow} y\right)\mathbb{P}_{\mathbf{p},q}\left(z\stackrel{A}{\longleftrightarrow} u\right)\\
&= \sum_{i=1}^N\sum_{\substack{\{y,z\}\in \\ \Delta S_x\cap E_i}}p_i\mathbb{P}_{\mathbf{p},q}\left(x\stackrel{S_x}{\longleftrightarrow} y\right)\sum_{u\in A\setminus S_x}\mathbb{P}_{\mathbf{p},q}\left(z\stackrel{A}{\longleftrightarrow} u\right)\\
&\leq \psi_{\mathbf{p},q}(x,S_x) \chi(A).
\end{split}
\]
We subsequently add the vertices in $S_x$ and use the trivial bound $\mathbb{P}_{\mathbf{p},q}\left(x\stackrel{A}{\longleftrightarrow}u\right)\leq 1$:
\[
\sum_{u\in A}\mathbb{P}_{\mathbf{p},q}\left(x\stackrel{A}{\longleftrightarrow}u\right)\leq |S_x|+\psi_{\mathbf{p},q}(x,S_x) \chi(A)
\]
The above inequality holds for any $x\in V$, so in particular it holds for the vertex which maximizes the left hand side. This vertex exists because the graph is quasi-transitive. We find
\begin{equation}
\max_{x\in V}\sum_{u\in A}\mathbb{P}_{\mathbf{p},q}\left(x\stackrel{A}{\longleftrightarrow}u\right)\leq \max_{x\in V}\big\{ |S_x|+\psi_{\mathbf{p},q}(x,S_x) \chi(A)\big\},
\end{equation}
so that
\begin{equation}
\chi(A)\leq \max_{x\in V}\big\{ |S_x|+\psi_{\mathbf{p},q}(x,S_x) \chi(A)\big\}.
\end{equation}
We conclude
\begin{equation}
\chi(A)\leq \frac{\max_{x\in V}|S_x|}{1-\max_{x\in V}\psi_{\mathbf{p},q}(x,S_x)},
\end{equation}
so that $\chi(A)$ is uniformly bounded from above. Therefore the final result follows by replacing $A$ with an exhausting sequence of subgraphs that tends to $V$.

\section{Supercritical Case: Existence of an Infinite Cluster}\label{sec:sup}
In this section we prove that $\Ppq(|\Ccal(x)| = \infty)>0$ when $q > q_\psi(\mathbf{p})$, for values of $\mathbf{p}\in \mathcal{A}$. If $q_\psi(\mathbf{p})=1$, then this statement is vacuous, so suppose $q_\psi(\mathbf{p})<1$ and let $q>q_\psi(\mathbf{p})$. We use the left-continuity of $\lambda \mapsto q_\psi(\mathbf{p}+\lambda \mathbf{1})$ implied by $\mathbf{p} \in \mathcal{A}$ to define an auxiliary point $(\hat{\mathbf{p}},\hat{q})$ in the $(\mathbf{p},q)$-space. Let $\varepsilon=\big(q-q_\psi(\mathbf{p})\big)/2$. Then there exists a $\delta>0$ such that
\[
q_\psi(\mathbf{p}+\lambda \mathbf{1})< q_\psi(\mathbf{p})+\varepsilon \qquad \textrm{for all}\quad \lambda \in [-\delta, 0].
\]
We define $\hat{\mathbf{p}}:=\mathbf{p}-\delta \mathbf{1}$ and $\hat{q}:= q_\psi(\hat{\mathbf{p}})$, so that $q-\hat{q}\geq \varepsilon$. Consider the line segment from $(\hat{\mathbf{p}},\hat{q})$ to $(\mathbf{p},q)$. We parametrise this line segment as
\[
\mathbf{r}(\lambda):=\begin{pmatrix}
\mathbf{p}+\lambda \mathbf{1}\\
q+\frac{\lambda}{\delta}(q-\hat{q})
\end{pmatrix},\qquad -\delta\leq\lambda\leq 0.
\]
Since $q_\psi(\hat{\mathbf{p}})=\hat{q}$ and since $\Gcol$ is quasi-transitive, there exists an $x\in V$ such that 
\[
\inf_{\substack{S\ssubset V: \\  S \owns x}}\psi_{\hat{\mathbf{p}},\hat{q}}(x,S)=1.
\]
Moreover, since $\psi_{\mathbf{p},q}(x,S)$ is increasing in $\mathbf{p}$ and $q$, we have for any $\lambda\in[-\delta,0]$ and for the same $x\in V$ that
\begin{equation}\label{eq:psi1}
\inf_{\substack{S\ssubset V: \\ S \owns x}}\psi_{\mathbf{r}(\lambda)}(x,S)\geq 1.
\end{equation}

Given a configuration $\omega \in \Omega$ and an edge $e$, we write $\omega_e$ for the element of $\Omega$ that satisfies $\omega_e(f) = \omega(f)$ for all $f \neq e$, and $\omega_e(e) = 1-\omega(e)$ (so $\omega_e$ is just $\omega$ with the status of the edge $e$ ``flipped''). Given an event $A$, we say that $e$ \emph{is pivotal for $A$ in $\omega$} if either $\omega \in A$ and $\omega_e \in A^c$ or $\omega \in A^c$ and $\omega_e \in A$. 
For two configurations $\omega, \omega' \in \Omega$ we write $\omega' \succcurlyeq \omega$ if $\omega'(e) \ge \omega(e)$ for all $e \in E$.
We say that an event $A$ is \emph{increasing} if for any $\omega \in A$ it holds that $\omega' \in A$ whenever $\omega' \succcurlyeq \omega$.

\emph{Russo's Formula}  \cite{Russo1981} states that the derivative of a percolation probability of an increasing event with respect to a parameter is equal to the expected number of pivotal edges associated with that parameter. More precisely, in our setting
\begin{lem}[Russo's Formula for $\Ppq$]\label{russo}
Let $A$ be an increasing event dependent only on edges in $\Lambda_n^v$ for some $n\in \mathbb{N}$ and $v\in V$, then for any $1\leq i\leq N$ we have
\[
\frac{\partial}{\partial p_i}\mathbb{P}_{\mathbf{p},q}(A)
  =\frac{1}{1-p_i}\sum_{e\in \Lambda_n^v\cap E_i}\mathbb{P}_{\mathbf{p},q}\big(\{e \text{ is pivotal for }A\}\cap A^c\big).
\]
\end{lem}
The proof is an easy adaptation of \cite[Lemma 3]{Russo1981}, so we leave it to the reader.

Fix $n\in \mathbb{N}$. We apply Russo's Formula to $x\longleftrightarrow \partial\Lambda_n^x$ (which is increasing), to find for all $1 \le i  \le N$,
\[
 \frac{\partial}{\partial p_i}\mathbb{P}_{\mathbf{p},q}(x\longleftrightarrow \partial\Lambda_n^x)=\frac{1}{1-p_i}\sum_{e\in \Lambda_n^x\cap E_i}\mathbb{P}_{\mathbf{p},q}(e \text{ pivotal}, x\nlongleftrightarrow \partial\Lambda_n^x),
\]
(where we abbreviated ``$e \text{ is pivotal for }x\longleftrightarrow \partial\Lambda_n^x$'' by ``$e$ pivotal'').
So we also have an expression for $\nabla \mathbb{P}_{\mathbf{p},q}(x\longleftrightarrow \partial\Lambda_n)$. We can integrate this gradient in the $(\mathbf{p},q)$-space along the straight line segment starting in $(\hat{\mathbf{p}},\hat{q})$ and ending in $(\mathbf{p},q)$. We use the Gradient Theorem to obtain
\begin{align}\label{eq:one}
\int_{-\delta}^0\nabla \mathbb{P}_{\mathbf{r}(\lambda)}(x\longleftrightarrow \partial\Lambda_n^x) \mathrm{d}\mathbf{r}(\lambda)=\mathbb{P}_{\mathbf{p},q}(x\longleftrightarrow \partial\Lambda_n^x)-\mathbb{P}_{\hat{\mathbf{p}},\hat{q}}(x\longleftrightarrow \partial\Lambda_n^x)\leq \mathbb{P}_{\mathbf{p},q}(x\longleftrightarrow \partial\Lambda_n^x).
\end{align}
On the other hand we can use the expression for $\nabla \mathbb{P}_{\mathbf{p},q}(x\longleftrightarrow \partial\Lambda_n^x)$ and the definition of the line integral along with the parametrization $\mathbf{r}(\lambda)$ to obtain
\begin{equation}\label{eq:other}
\begin{split}
\int_{-\delta}^0 \nabla \mathbb{P}_{\mathbf{r}(\lambda)}(x\longleftrightarrow & \partial\Lambda_n^x) \mathrm{d}\mathbf{r}(\lambda)=\int_{-\delta}^0 \nabla \mathbb{P}_{\mathbf{r}(\lambda)}(x\longleftrightarrow \partial\Lambda_n^x)\cdot \mathbf{r}'(\lambda)\,\mathrm{d}\lambda\\
&\geq\sum_{i=1}^{N}\int_{-\delta}^0 \left(\frac{\varepsilon}{\delta}\wedge 1 \right) \sum_{e\in \Lambda_n^x\cap E_i}\mathbb{P}_{\mathbf{r}(\lambda)}(e \text{ pivotal}, x\nlongleftrightarrow \partial\Lambda_n^x)\,\mathrm{d}\lambda,
\end{split}
\end{equation}
since $\mathbf{r}'(\lambda)=\big(\mathbf{1},\frac{1}{\delta}(q-\hat{q})\big)^T$. We define $c:=\frac{\varepsilon}{\delta}\wedge 1$, so that $c>0$.

We now define the random subset of $\Lambda_n^x$
\begin{equation}
\mathscr{S}:=\big\{y\in \Lambda_n^x \: : \: y\nlongleftrightarrow \partial\Lambda_n^x\big\}.
\end{equation}
The boundary of $\mathscr{S}$ are the vertices of $\Lambda_n^x$ for which all neighbours that are not in $\mathscr{S}$ are connected to $\partial\Lambda_n^x$. If $x\nlongleftrightarrow \partial\Lambda_n^x$, then $x\in \mathscr{S}$, so if we sum over all possible values of $\mathscr{S}$ we find
\begin{equation}\label{eq:combined}
\mathbb{P}_{\mathbf{p},q}(x\longleftrightarrow \partial\Lambda_n^x)\geq c\sum_{i=1}^{N}\int_{-\delta}^0 \sum_{\substack{S\subset \Lambda_n^x:\\  S \owns x}}\sum_{\substack{e\in  \Lambda_n^x\cap E_i}}\mathbb{P}_{\mathbf{r}(\lambda)}(e \text{ pivotal}, \mathscr{S}=S)\,\mathrm{d}\lambda
\end{equation}
When we know $\mathscr{S}=S$ we know that the pivotal edges for the event $\{x\longleftrightarrow\partial\Lambda_n^x\}$ are exactly the edges $\{y,z\}$ on the edge-boundary $\Delta S$ of $S$ that are connected to $x$, i.e., $y\in S$, $z\not\in S$ and $x\longleftrightarrow y$. We can sum over these edges to obtain
\[
\mathbb{P}_{\mathbf{p},q}(x\longleftrightarrow \partial\Lambda_n^x)\geq c\sum_{i=1}^{N}\int_{-\delta}^0 \sum_{\substack{S\subset \Lambda_n^x:\\ S \owns x}}\:\sum_{\substack{\{y,z\}\in \\ \Delta S\cap E_i}}\mathbb{P}_{\mathbf{r}(\lambda)}\left(x\stackrel{S}{\longleftrightarrow} y, \mathscr{S}=S\right)\mathrm{d}\lambda.
\]

The occurrence of the event $\{\mathscr{S}=S\}$ can be determined from the configuration of the edges outside of $S$. This can be done by exploring from the boundary of $\Lambda_n^x$: Let $\mathscr{T}$ be the set of vertices in $\Lambda_n^x$ that are connected to $\partial\Lambda_n^x$ using only edges in $\Lambda_n^x\backslash S$. Then $\{\mathscr{S}=S\}=\{\mathscr{T}=\Lambda_n^x\backslash S\}$. We conclude that the event $\{\mathscr{S}=S\}$ is determined by the configuration of the edges outside $S$, and is therefore independent of $\{x\stackrel{S}{\longleftrightarrow} y\}$. We find
\begin{equation}\label{eq:moeilijk}
\mathbb{P}_{\mathbf{p},q}(x\longleftrightarrow \partial\Lambda_n^x)\geq c\sum_{i=1}^{N}\int_{-\delta}^0 \sum_{\substack{S\subset \Lambda_n^x:\\  S \owns x}}\:\sum_{\substack{\{y,z\}\in \\ \Delta S\cap E_i}}\mathbb{P}_{\mathbf{r}(\lambda)}\left(x\stackrel{S}{\longleftrightarrow} y\right)\mathbb{P}_{\mathbf{r}(\lambda)}\left(\mathscr{S}=S\right)\mathrm{d}\lambda.
\end{equation}
By (\ref{eq:psi1}) we have that $\psi_{\mathbf{r}(\lambda)}(x,S_x)\geq 1$ for all $x\in S_x\subset \mathbb{Z}^d$. Using the definition of $\psi_{\mathbf{p},q}(x,S)$ we obtain
\begin{align}
\mathbb{P}_{\mathbf{p},q}(x\longleftrightarrow \partial\Lambda_n^x)&\geq c \int_{-\delta}^0 \sum_{\substack{S\subset \Lambda_n^x:\\  S \owns x}}\psi_{\mathbf{r}(\lambda)}(x,S_x)\mathbb{P}_{\mathbf{r}(\lambda)}(\mathscr{S}=S)\,\mathrm{d}\lambda\nonumber\\
&\geq c \int_{-\delta}^0\mathbb{P}_{\mathbf{r}(\lambda)}(x\nlongleftrightarrow \partial\Lambda_n^x)\,\mathrm{d}\lambda\nonumber\\
&\geq c\delta \mathbb{P}_{\mathbf{p},q}(x\nlongleftrightarrow \partial\Lambda_n^x).
\end{align}
Finally we obtain
\begin{align}
\mathbb{P}_{\mathbf{p},q}(x\longleftrightarrow \partial\Lambda_n^x)\geq \frac{\varepsilon \wedge \delta}{1+\varepsilon \wedge \delta }>0.
\end{align}
The way we chose $\varepsilon$ and $\delta$ was independent of $n$, so letting $n$ tend to infinity, it follows that $\Ppq(|\Ccal(x)| = \infty) >0$. Now let $y\in V$, since $G$ is connected it holds that $\mathbb{P}_{\mathbf{p},q}(x\longleftrightarrow y)>0$. It follows that $\Ppq(|\Ccal(y)| = \infty) >0$ for all $y\in V$ and all $q > q_\psi(\mathbf{p})$, as claimed.

\section{Subcritical case: exponential decay}\label{sec:exp}
In this section we prove that for any $\mathbf{p}\in[0,1]^{N-1}$, \eqref{eq:thma} holds if $q<q_\psi(\mathbf{p})$. If $q<q_\psi(\mathbf{p})$, it follows that for any $x\in V$ there exists a finite subset $S_x\ssubset V$ with $x\in S_x$ such that $\psi_{\mathbf{p},q}(x,S_x)<1$. Furthermore, since $\Gcol$ is quasi transitive there are only finitely many different types of vertices, so we can find an $L\in \mathbb{N}$ such that $S_x\subset\Lambda_L^x$ for any $x\in V$. Now let $x\in V$ be given and fix a set $S$ with $x\in S$ and $\psi_{\mathbf{p},q}(x,S)<1$. Since $\Gcol$ is quasi-transitive there exists a $u\in V$ such that
\[
\min_{\substack{S\subset \Lambda_L^u: \\ S \owns u}}\psi_{\mathbf{p},q}(u,S)=\sup_{v\in V}\min_{\substack{S\subset \Lambda_L^v: \\ S \owns v}}\psi_{\mathbf{p},q}(v,S).
\]
The above step is (one of the places) where the proof would fail if the graph was not quasi-transitive. Because of the quasi-transitivity the above supremum is really a maximum that is attained for some $u\in V$. Without quasi-transitivity it is not clear if the supremum would be attained in some vertex, or at infinity.

Since $q<q_\psi(\mathbf{p})$ there must exist $\varepsilon>0$ such that
\[
\min_{\substack{ S\subset  \Lambda_L^u: \\ S \owns u}}\psi_{\mathbf{p},q}(u,S)=1-\varepsilon.
\]
Define the random set
\[
\mathscr{C}_S:=\{z\in S: x\stackrel{S}{\longleftrightarrow} z\}.
\]
Now if we let $k\in \mathbb{N}$ and we suppose $x\longleftrightarrow \partial \Lambda^x_{kL}$, then we know that there exists an edge $\{y,z\}\in \Delta S$ such that $x \stackrel{S}{\longleftrightarrow} y$, $\{y,z\}$ is open and $z\stackrel{\mathscr{C}_S^c}{\longleftrightarrow} \partial\Lambda^x_{kL}$. So by summing over all possible edges in $\Delta S$ and over all possible values of $\mathscr{C}_S$ we obtain
\[
\mathbb{P}_{\mathbf{p},q} (x\longleftrightarrow \partial \Lambda^x_{kL})\leq\sum_{C\subset S}\sum_{\{y,z\}\in \Delta S}\mathbb{P}_{\mathbf{p},q}\left(\left\{x \stackrel{S}{\longleftrightarrow} y,\mathscr{C}_S=C\right\}, \{y,z\}\text{ open, }z\stackrel{C^c}{\longleftrightarrow} \partial\Lambda^x_{kL}\right).
\]
The three events in the above inequality depend on disjoint sets of edges: $\big\{x \stackrel{S}{\longleftrightarrow} y,\mathscr{C}_S=C\big\}$ only depends on edges with both end points in $C$, the event $z\stackrel{C^c}{\longleftrightarrow} \partial\Lambda^x_{kL}$ depends on edges with both end points in $C^c$ and the edge $\{y,z\}$ has neither both end points in $C$ nor in $C^c$. So these events are independent and we obtain
\begin{align}\label{eq:prf1}
	\mathbb{P}_{\mathbf{p},q}(x\longleftrightarrow \partial \Lambda^x_{kL})\leq \sum_{i=1}^{N} p_i\sum_{C\subset S}\sum_{\substack{\{y,z\}\in \\ \Delta S\cap E_i}}\Ppq \left(x \stackrel{S}{\longleftrightarrow} y,\mathscr{C}_S =C\right)\Ppq \left(z\longleftrightarrow \partial\Lambda^x_{kL}\right).
\end{align}
Since $z\in \Lambda^x_L$ we have $\{z \longleftrightarrow \partial \Lambda_{(k-1)L}^z\} \subset \{z\longleftrightarrow \partial\Lambda^x_{kL}\}$, so we find
\[
\mathbb{P}_{\mathbf{p},q}\left(z\longleftrightarrow \partial\Lambda^x_{kL}\right)\leq \mathbb{P}_{\mathbf{p},q}\left(z\longleftrightarrow \partial\Lambda^z_{(k-1)L}\right).
\]
Furthermore since $\Gcol$ is quasi-transitive there exists a $w\in V$ such that
\[
\mathbb{P}_{\mathbf{p},q}\left(w\longleftrightarrow \partial\Lambda^w_{(k-1)L}\right)=\sup_{v \in V} \mathbb{P}_{\mathbf{p},q}\left(v\longleftrightarrow \partial\Lambda^v_{(k-1)L}\right)\geq \mathbb{P}_{\mathbf{p},q}\left(z\longleftrightarrow \partial\Lambda^z_{(k-1)L}\right).
\]
So that for this $w\in V$
\[
\mathbb{P}_{\mathbf{p},q}\left(z\longleftrightarrow \partial\Lambda^x_{kL}\right)\leq\mathbb{P}_{\mathbf{p},q}\left(w\longleftrightarrow \partial\Lambda^w_{(k-1)L}\right).
\]
This bound can now be used in (\ref{eq:prf1}) to obtain
\begin{equation}\label{eq:bnd}
\begin{split}
\mathbb{P}_{\mathbf{p},q}(x\longleftrightarrow \partial \Lambda^x_{kL})&\leq \sum_{i=1}^{N} p_i\sum_{C\subset S}\sum_{\substack{\{y,z\}\in \\ \Delta S\cap E_i}}\Ppq \left(x \stackrel{S}{\longleftrightarrow} y,\mathscr{C}_S=C\right)\mathbb{P}_{\mathbf{p},q}\left(w\longleftrightarrow \partial\Lambda^w_{(k-1)L}\right)\\
&\leq \psi_{\mathbf{p},q}(x,S)\Ppq (w\longleftrightarrow \partial\Lambda^w_{(k-1)L})\\
&\leq (1-\varepsilon)\mathbb{P}_{p,q}(w\longleftrightarrow \partial\Lambda^w_{(k-1)L}),
\end{split}
\end{equation}
where, in the second step we used the definition of $\psi_{\mathbf{p},q}(x,S)$, \eqref{e:psidef}.
Iteration of the above inequality now gives the desired exponential decay:
\begin{equation}
\mathbb{P}_{\mathbf{p},q}(x\longleftrightarrow \partial \Lambda^x_{kL})\leq (1-\varepsilon)^k=\exp\big(\log(1-\varepsilon)k\big)=\exp(-ck),
\end{equation}
where we set $c :=  \log(1-\varepsilon)$.

\section*{Acknowledgments}

TH is supported by the Netherlands Organisation for Scientific Research 
(NWO) through Gravitation-grant NETWORKS-024.002.003.

\bibliographystyle{abbrv}
\bibliography{ip}

\end{document}